\documentclass[journal]{IEEEtran}
\usepackage{amsmath,amssymb,mathrsfs,amsfonts,graphicx,latexsym,exscale,cmmib57,dsfont,bbm,amscd,ulem}
\usepackage{mathrsfs}

\newcommand{\bP}{\boldsymbol{P}}
\newcommand{\bp}{\boldsymbol{p}}
\newcommand{\bS}{\boldsymbol{S}}
\newcommand{\bxi}{\boldsymbol\xi}
\newcommand{\bcdot}{\boldsymbol\cdot}
\newenvironment{proof}{{\bf Proof:}}{\hfill$\Box$\medskip}
\newtheorem{thm}{\textbf{Theorem}}[section]

\newtheorem{lem}[thm]{\textbf{Lemma}}
\newtheorem{prop}[thm]{\textbf{Proposition}}

\newtheorem{defn}[thm]{\textbf{Definition}}

\newtheorem{remark}[thm]{\textbf{Remark}}

\newcommand{\EQ}{\begin{equation}\begin{array}{lllllllll}}
\newcommand{\EE}{\end{array}\end{equation}}
\newcommand{\MT}{\left[ \begin{array}{ccccccccc}}
\newcommand{\EM}{\end{array}\right]}

\newcommand{\eq}{\begin{equation}\begin{array}{lclllllllllllllll}}
\newcommand{\ee}{\end{array}\end{equation}}
\newcommand{\bmt}{\left[ \begin{array}{ccccccccc}}
\newcommand{\emt}{\end{array}\right]}
\newcommand{\bea}{\begin{eqnarray}}
\newcommand{\eea}{\end{eqnarray}}
\newcommand{\bean}{\begin{eqnarray*}}
\newcommand{\eean}{\end{eqnarray*}}

\ifCLASSINFOpdf
\else
\fi
\begin{document}
%
\title{Pointwise Stabilization of Discrete-time Stationary Matrix-valued Markovian Chains}
%
%
%

\author{Xiongping~Dai, Yu~Huang,
        and~Mingqing~Xiao,~\IEEEmembership{Senior Member,~IEEE}
\thanks{X. Dai is with Nanjing University, Nanjing 210093, People's Republic of China.  Email: xpdai@nju.edu.cn.}
\thanks{Y. Huang is with Zhongshan (Sun Yat-Sen) University, Guangzhou 510275, People's Republic of China. Email: stshyu@mail.sysu.edu.cn.}
\thanks{M. Xiao is with Southern Illinois University, Carbondale, IL 62901, USA. Email: mxiao@siu.edu.}}
\maketitle

\begin{abstract}
We study the pointwise stabilizability of a discrete-time, time-homogeneous, and stationary Markovian jump linear system. By using measure theory, ergodic theory and a splitting theorem of state space we show in a relatively simple way that if the system is essentially product-bounded, then it is pointwise convergent if and only if it is pointwise exponentially convergent, which provides an important characteristic of pointwise convergence under the framework of symbolic dynamics.
\end{abstract}

\begin{IEEEkeywords}
Discrete-time Markovian jump systems; pointwise convergence; pointwise exponential convergence.
\end{IEEEkeywords}

%
\IEEEpeerreviewmaketitle

\section{Introduction}
%
%
%
%
\IEEEPARstart{I}{n} this technical note, we will study the pointwise stabilization of a discrete-time, matrix-valued, and stationary Markovian chain. Let $\bS=\{S_1,\dotsc,S_K\}$ be an arbitrary set that consists of $K$ real $d\times d$ matrices, where $K$ and $d$ both are integers with $2\le K<\infty$ and $2\le d<\infty$. Set $\mathbb{K}=\{1,\dotsc,K\}$ equipped with the discrete topology. The system $\bS$ is said to be:
\begin{itemize}
\item ``pointwise convergent" if for each $x\in\mathbb{R}^{1\times d}$, there is an infinite switching sequence, say $i_{\bcdot}(x)\colon \mathbb{N}\rightarrow\mathbb{K}$ such that
    \bean
    \lim_{n\to+\infty}xS_{i_1(x)}\dotsm S_{i_n(x)}=\mathbf{0};
    \eean
\item ``pointwise exponentially convergent" if for each initial state $x\in\mathbb{R}^{1\times d}$, there is an infinite switching sequence, say $i_{\bcdot}(x)\colon \mathbb{N}\rightarrow\mathbb{K}$ such that
    \bean
    \limsup_{n\to+\infty}\frac{1}{n}\log\|xS_{i_1(x)}\dotsm S_{i_n(x)}\|<0.
    \eean
\end{itemize}

Here and in the sequel, $\mathbb{N}=\{1,2,\dotsc\}$, $\mathbf{0}$ stands for the origin of $\mathbb{R}^{d_1\times d_2}$, and by $\|\cdot\|$ we denote the usual Euclidean norm on $\mathbb{R}^{1\times d}$ defined by
\bean
\|x\|=\sqrt{x x^\mathrm{T}}\quad\forall x\in\mathbb{R}^{1\times d}
\eean
and also the matrix norm on $\mathbb{R}^{d\times d}$ compatible with the $d$-dimensional row-vector norm $\|\cdot\|$ on $\mathbb{R}^{1\times d}$.

It is clear that the notion of ``convergence" here is abused as it is referring to the usual approach to ``convergence" in the stability theory that requires convergence for any (or almost any) switching sequences $i_{\bcdot}=(i_{n})_{n=1}^{+\infty}\in\mathbb{K}^\mathbb{N}$.
Here our convergence takes place only for some desired index sequence $i_{\bcdot}(x)$. Further, $\bS$ is said to be
\begin{itemize}
\item ``consistently convergent" (also called ``uniformly convergent" in, for example, \cite{SU94} and \cite{CRS99}), if the switching sequence $i_{\bcdot}(x)$ in the pointwise convergence can be taken independent of the initial state $x$; that is to say, there exists a switching sequence $i_{\bcdot}\colon\mathbb{N}\rightarrow\mathbb{K}$ such that
    \bean
    \lim_{n\to+\infty}S_{i_1}\dotsm S_{i_n}=\mathbf{0}.
    \eean
\end{itemize}

These concepts arise and have been  studied naturally in the theory of multi-rate sampled-data control systems and multi-modal linear control systems and for some control optimization problems in, for example, see \cite{Stan79, BCS88, SU94, GR97, CRS99, Sun04, Sun06, LD06, LD07, Dai-JDE}. In general, pointwise convergence does not imply the pointwise exponential convergence.

In this paper, we consider the random version of the above important concepts driven by discrete-time stationary matrix-valued Markovian chains. Under the framework of symbolic dynamics, we can prove our main result in a relatively simple way by using a dichotomy theorem of the state space $\mathbb{R}^{1\times d}$.

Let $(\Omega,\mathscr{F},\mathbb{P})$ be a probability space throughout the sequel of this note and let
\bean
\bxi=(\xi_n)_{n=1}^{+\infty},\quad \textrm{where }\xi_n\colon\Omega\rightarrow\mathbb{K},
\eean
be a discrete-time, time-homogeneous, and stationary $(\bp,\bP)$-Markovian chain, which naturally induces a ``Markovian jump linear system" based on $\bS$ as follows:\footnote{In almost all available literature, peoples directly consider the matrix-valued Markovian process
\bean
\bxi=(\xi_n)_{n\ge1},\quad \textrm{where } \xi_n\colon\Omega\rightarrow\bS\textrm{ instead of }\xi_n\colon\Omega\rightarrow\mathbb{K},
\eean
as formulated in Abstract. This looks more concise. However, our treatment presented here enables us to employ the abstract theory of symbolic dynamical systems. See Section~\ref{sec2} below.\label{foot1}}
\bean
x_{n}=x{S}_{\xi_1(\omega)}\cdots{S}_{\xi_n(\omega)},\quad x\in\mathbb{R}^{1\times d}, n\ge1,\; \omega\in\Omega.
\eean
Here $\bp=(p_1,\dotsc,p_K)\in\mathbb{R}^{1\times K}$ is the initial distribution of $\bxi$ and $\bP=[p_{ij}]\in\mathbb{R}^{K\times K}$ is the Markov transition probability matrix of $\bxi$ satisfying the stationary distribution property:
\bean
\bp\bP=\bp.
\eean
For any sample $\omega\in\Omega$, there corresponds an infinite switching sequence
\bean
\bxi(\omega)\colon \mathbb{N}\rightarrow\mathbb{K};\; n\mapsto\xi_n(\omega)
\eean
that is named a ``trajectory" of the Markovian chain $\bxi$ in the textbooks of stochastic processes. Moreover, for a given switching sequence $i_{\bcdot}=(i_n)_{n=1}^{+\infty}\in\mathbb{K}^\mathbb{N}$, it is not necessary to have some sample $\omega\in\Omega$ satisfying $\xi_n(\omega)=i_n$ for all $n\ge1$. However, under the framework of probability, we are only interested in the events of positive-probability of $\mathbb{P}$. Thus, motivated by this consideration, we introduce the following concepts.

\begin{defn}\label{def1.1}
The Markovian jump linear system $(\bS,\bxi)$ is said to be:
\begin{enumerate}
\item[$\mathrm{(a)}$] ``pointwise convergent", if to any initial state $x\in\mathbb{R}^{1\times d}$, there corresponds an $\Omega_{x}\in\mathscr{F}$ with $\mathbb{P}(\Omega_{x})>0$ such that
\bean
    x{S}_{\xi_1(\omega)}\dotsm{S}_{\xi_n(\omega)}\to\mathbf{0}\;\textrm{ as }\;n\to+\infty,\quad \forall \omega\in\Omega_{x};
\eean

\item[$\mathrm{(b)}$] ``pointwise exponentially convergent", if to any initial state $x\in\mathbb{R}^{1\times d}$, there corresponds an $\Omega_{x}^\prime\in\mathscr{F}$ with $\mathbb{P}(\Omega_{x}^\prime)>0$ such that
\bean
    \limsup_{n\to+\infty}\frac{1}{n}\log\|x{S}_{\xi_1(\omega)}\dotsm{S}_{\xi_n(\omega)}\|<0,\quad \forall \omega\in\Omega_{x}^\prime;
\eean

\item[$\mathrm{(c)}$] ``consistently exponentially convergent", if there exists a subset $\Omega^{\prime\prime}\in\mathscr{F}$ with $\mathbb{P}(\Omega^{\prime\prime})>0$ such that
\bean
    \limsup_{n\to+\infty}\frac{1}{n}\log\|{S}_{\xi_1(\omega)}\dotsm{S}_{\xi_n(\omega)}\|<0,\quad \forall \omega\in\Omega^{\prime\prime}.
\eean
\end{enumerate}
\end{defn}

Here ``consistently" only means that the choice of the driving sample $\omega\in\Omega$ is independent of any initial state $x\in\mathbb{R}^{1\times d}$. We should notice that the consistent exponential convergence of $(\bS,\bxi)$ over $\omega$ is essentially weaker than the so-called ``uniform exponential stability" over $\omega$ which means that there are constants $C>0$ and $0<\gamma<1$ such that
\bean
\|{S}_{\xi_{\ell+1}(\omega)}\dotsm{S}_{\xi_{\ell+m}(\omega)}\|\le C\gamma^m\quad\forall\ell\ge0\textrm{ and }m\ge1.
\eean
In addition, since $\bxi$ does not need to be irreducible (or equivalently, not need to be ergodic; see Section~\ref{sec2}), it may not be realistic to require $\mathbb{P}(\Omega^{\prime\prime})=1$ in general.
It is obvious that (c)$\Rightarrow$(b)$\Rightarrow$(a), but not vice versa in general. For example, let $\bS$ consist of
\bean
S_1=\left[\begin{matrix}\frac{1}{2}&0\\0&2\end{matrix}\right],\; S_2=\left[\begin{matrix}\frac{\sqrt{3}}{2}&\frac{1}{2}\\-\frac{1}{2}&\frac{\sqrt{3}}{2}\end{matrix}\right]
\eean
and
\bean
\xi_n(\omega)\in\{1,2\}\quad\textrm{ with }\bp=({1}/{2},{1}/{2}),\;\bP=\left[\begin{matrix}\frac{1}{2}&\frac{1}{2}\\\frac{1}{2}&\frac{1}{2}\end{matrix}\right];
\eean
then $(\bS,\bxi)$ is pointwise convergent but not consistently exponentially convergent.

\begin{remark}\label{rem1.2}
Because considering a deterministic sample $\omega$ in $\Omega$ makes no sense in probability theory and random/stochastic stability theory, it is necessary to require the property of positive probability: $\mathbb{P}(\Omega_{x})>0$, $\mathbb{P}(\Omega_{x}^\prime)>0$, and $\mathbb{P}(\Omega^{\prime\prime})>0$, in Definition~\ref{def1.1}. That means all $\Omega_{x}, \Omega_{x}^\prime$ and $\Omega^{\prime\prime}$ to be non-negligible events that characterize the convergent properties under the framework of probability.
\end{remark}

\begin{remark}\label{rem1.3}
If $(\bS,\bxi)$ is ergodic $($in fact, we can always assume this by an ergodic decomposition Theorem~\ref{thm2.3} below$)$, then it may be required that $\mathbb{P}(\Omega^{\prime\prime})=1$ in Definition~\ref{def1.1}.$\mathrm{(c)}$; this is because from the Subadditive Ergodic Theorem, see e.g. \cite[Theorem~10.1]{Wal82}, the stable set of $(\bS,\bxi)$
\bean
\Omega_{\textrm{stable}}=\left\{\omega\colon \lim_{n\to+\infty}\frac{1}{n}\log\|{S}_{\xi_1(\omega)}\dotsm{S}_{\xi_n(\omega)}\|\equiv\chi<0\right\},
\eean
which is invariant, has the $\mathbb{P}$-measure either $1$ or $0$. However, for Definition~\ref{def1.1}.$\mathrm{(a)}$, in general $\Omega_x$ is not necessarily invariant and hence
\bean
\mathbb{P}(\Omega_{x})>0\not\Rightarrow\mathbb{P}\left\{\omega\colon x{S}_{\xi_1(\omega)}\dotsm{S}_{\xi_n(\omega)}\to\mathbf{0}\right\}=1.
\eean
This is just one of the essential challenges for $\mathrm{(a)}\not\Rightarrow \mathrm{(c)}$. For more details, see Section~\ref{sec3} below.
\end{remark}

By using measure theory, ergodic theory and a splitting theorem of the state space $\mathbb{R}^{1\times d}$, we will mainly show in Section~\ref{sec3} that \textit{if $(\bS,\bxi)$ is essentially nonuniformly product bounded then $(\bS,\bxi)$ is pointwise convergent if and only if it is pointwise exponentially convergent}. See Theorem~\ref{thm3.1} below for an equivalent formulation in terms of symbolic dynamics.

We notice that if $\bS$ is irreducible\footnote{The matrix family $\bS$ is said to be ``irreducible" if there are no common, proper, nonempty, and invariant subspaces of $\mathbb{R}^{1\times d}$, for each member of $\bS$. See, e.g., \cite{Bar88}.\label{foot2}}
with the joint spectral radius $\hat{\rho}(\bS)=1$, then $(\bS,\bxi)$ is uniformly product bounded and hence essentially nonuniformly product bounded; see, e.g., N.~Barabanov~\cite{Bar88} also X.~Dai~\cite{Dai-JMAA}. The boundedness condition, also named as ``Lyapunov stability" in ODE, is both practically important and academically challenging \cite{MS, BL00, Hart02, BL05, Sun08}. Indeed, it is desirable in many practical issues and is closely related to periodic solutions and limit cycles; see, e.g., \cite{Bar88-SMJ} and \cite{Bar93}.

Such an equivalence theorem will play a key role in creating upper bounds, finding convergence rates and exploiting other basic system properties for Markovian jump linear systems, as done in the deterministic case, for example, in \cite{Sun06} and \cite{HC08}.

To prove our main result (Theorem~\ref{thm3.1}), we need to make use of two important tools. One is the ergodic theory of Markovian chains established in Section~\ref{sec2},
and the other is a dichotomy decomposition theorem (lemma~\ref{lem3.4}) proved in \cite{DHX-pre}. We will end this paper with concluding remarks in Section~\ref{sec4}.
\section{Ergodic theory of Markovian chains}\label{sec2}%
In this section, we will introduce the framework of the stationary Markovian chains in ergodic theory, which will be used in the proof of our main
result.

Let $\bxi=(\xi_n)_{n=1}^{+\infty}$ where $\xi_n\colon\Omega\rightarrow\mathbb{K}$, be a discrete-time, \textit{time-homogeneous}, and \textit{stationary} $(\bp,\bP)$-Markovian chain, defined on the probability space $(\Omega,\mathscr{F},\mathbb{P})$ with the finite state-space
$\mathbb{K}=\{1,\dotsc,K\}$
that is equipped with the discrete topology. Notice here that ``time-homogeneity" means that the transition probabilities,
\bean
\mathbb{P}[\xi_{n+1}=j\,\pmb{|}\,\xi_n=i]=p_{ij}\; \forall i,j\in\mathbb{K},
\eean
where $\bP=[p_{ij}]\in\mathbb{R}^{K\times K}$, all do not depend upon the time $n$; and ``stationary" means
\bean
\bp\bP=\bp.
\eean
This implies that $\xi_1,\xi_2,\dotsc$ are identically distributed random variables. However, they are not necessarily independent. It is straightforward to see
\begin{equation}\label{eq2.1}
\mathbb{P}[\xi_1=i_1, \dotsc, \xi_n=i_n]=p_{i_1}p_{i_1i_2}\dotsm p_{i_{n-1}i_n}
\end{equation}
for any words
$(i_1,\dotsc,i_n)\in\mathbb{K}^n$
of finite-length $n\ge2$, from the Chapman-Kolmogorov equation. Here the probability row-vector
$\bp=(p_1,\dotsc,p_K)\in\mathbb{R}^{1\times K}$
is the initial distribution of the Markovian chain $\bxi$, i.e.,
$\mathbb{P}[\xi_1=k]=p_k$
for all $1\le k\le K$. Hereafter, assume $\bp>0$, i.e., $p_k>0\; \forall k\in\mathbb{K}$; otherwise, we only need to replace the state-space $\mathbb{K}$ of the Markovian chain $\bxi$ with $\mathbb{K}\setminus\{k\,|\,p_{k}=0\}$ by the standard treatment.

The transition probability matrix $\bP$ of $\bxi$ is called ``irreducible" if for any pair of states $i,j\in\mathbb{K}$, there is some $n=n(i,j)\ge1$ such that $p_{ij}^{(n)}>0$, where $p_{ij}^{(n)}$ is the $(i,j)$-coordinate element of the $n$-time product matrix $\bP^n$. It is worth to mention here that although this ``irreducibility" of $\bP$ has nothing in conceptual common with the irreducibility explained to a family of matrices $\bS$ in Footnote~\ref{foot2} before, but in fact the irreducibility of Markovian chains does have a very close connection to that as defined in Footnote~\ref{foot2}; see, for example, \cite{Sen}.

We denote by $\varSigma_K^+$ the set of all infinite switching sequences $i_{\bcdot}\colon\mathbb{N}\rightarrow\mathbb{K};\; n\mapsto i_n$. Here it is convenient to place the variables $1,2,\dotsc$ at the subscript position. Then, under the infinite product topology that can be generated by the cylinders
\bean
[i_1^\prime,\dotsc,i_\ell^\prime]=\left\{i_{\bcdot}\in\varSigma_K^+\,|\,i_1=i_1^\prime, \dotsc, i_\ell=i_\ell^\prime\right\}
\eean
for all $\ell\ge1$ and any words $(i_1^\prime,\dotsc,i_\ell^\prime)\in\mathbb{K}^\ell$ of finite-length $\ell$, $\varSigma_K^+$ is a compact topological space as well as
the one-sided Markovian shift transformation
\bean
\theta\colon\varSigma_K^+\rightarrow\varSigma_K^+;\quad i_{\bcdot}=(i_n)_{n=1}^{+\infty}\mapsto i_{\bcdot+1}=(i_{n+1})_{n=1}^{+\infty}
\eean
is continuous and surjective. But it is not injective, which presents some challenges since the stable manifold and/or the center manifold may not be invariant for the state trajectories (see detailed discussion in lemma~\ref{lem3.4}).

By the joint random variable setting
\bean
\varXi\colon\Omega\rightarrow\varSigma_K^+;\quad \omega\mapsto\bxi(\omega)=(\xi_n(\omega))_{n=1}^{+\infty},
\eean
we can obtain a natural probability distribution, called the ``$(\bp,\bP)$-Markovian measure" and simply write as $\mu_{\bp,\bP}$, on $\varSigma_K^+$, which is such that for any $n\ge1$,
\begin{equation}\label{eq2.2}
\mu_{\bp,\bP}([i_1,\dotsc,i_n])=\begin{cases}
p_{i_1}& \textrm{if }n=1;\\
p_{i_1}p_{i_1i_2}\dotsm p_{i_{n-1}i_n}& \textrm{if }n\ge2,
\end{cases}
\end{equation}
for all cylinder sets $[i_1,\dotsc,i_n]\subset\varSigma_K^+$. So, $\mu_{\bp,\bP}=\mathbb{P}\circ\varXi^{-1}$ from (\ref{eq2.1}) and (\ref{eq2.2}).

It should be noted here that $\mu_{\bp,\bP}$ is not necessarily equal to the infinite product of the initial distribution $\bp$ of $\bxi$, for $\xi_1,\xi_2,\dotsc$ need not be independent each others. In addition, $\varXi(\Omega)\not=\varSigma_K^+$ in general, unless $p_{ij}>0$ for all $i,j\in\mathbb{K}$.

The following is a known result that will be used for our arguments later.

\begin{lem}[P.~Walters~\cite{Wal82}]\label{lem2.1}
Let $\bxi=(\xi_n)_{n=1}^{+\infty}\colon\Omega\rightarrow\varSigma_K^+$ be a discrete-time, time-homogeneous, and stationary $(\bp,\bP)$-Markovian chain. Then,
\begin{enumerate}
\item the Markovian shift transformation $\theta\colon\varSigma_K^+\rightarrow\varSigma_K^+$ preserves the $(\bp,\bP)$-Markovian measure $\mu_{\bp,\bP}$; that is,
\bean
\mu_{\bp,\bP}(B)=\mu_{\bp,\bP}\circ\theta^{-1}(B)\quad\forall B\in\mathscr{B}_{\varSigma_K^+};
\eean
\item $\bP$ is irreducible if and only if $\mu_{\bp,\bP}$ is $\theta$-ergodic on $\varSigma_K^+$; that is, for all Borel subsets $B\subset\varSigma_K^+$, the equality $\mu\left((\theta^{-1}(B)\setminus B)\cup(B\setminus\theta^{-1}(B))\right)=0$
implies that $\mu(B)=1$ or $0$.
\end{enumerate}
\end{lem}
Here $\mathscr{B}_{\varSigma_K^+}$ is the standard Borel $\sigma$-field of the space $\varSigma_K^+$.

Since the Markovian transition probability matrix $\bP$ is not necessarily irreducible in our situation, we need to consider the ergodic decomposition of the $(\bp,\bP)$-Markovian probability $\mu_{\bp,\bP}$.  A state $k\in\mathbb{K}$ is called ``recurrent" for $\bxi$, if the conditional probability
\bean
\mathbb{P}[\omega\in\Omega\colon\exists n_\ell\nearrow+\infty\textrm{ s.t. }\xi_{n_\ell}(\omega)=k\,\pmb{|}\,\xi_1=k]=1.
\eean
If $k\in\mathbb{K}$ is not recurrent, then it is called ``non-recurrent" or ``transient". Any two states $k_1,k_2$, each accessible to the other, i.e, $p_{k_1k_2}^{(m)}>0$ and $p_{k_2k_1}^{(n)}>0$ for some pair $m,n\ge1$, are said to be ``communicative" and we write $k_1\leftrightsquigarrow k_2$. The concept of $\leftrightsquigarrow$ is an equivalence relationship.

Then according to the classical theory of stochastic processes, for example, \cite{Ch60}, there exists the following basic partition of the states:
\bean
\mathbb{K}=\mathbb{K}_0\cup\mathbb{K}_1\cup\cdots\cup\mathbb{K}_r
\eean
such that:
\begin{itemize}
\item $\mathbb{K}_0$ consists of all the non-recurrent states of the Markovian chain $\bxi$;

\item each $\mathbb{K}_i, 1\le i\le r$, is closed and communicative, i.e., for any $k,k^\prime\in \mathbb{K}_i$ and $k^{\prime\prime}\not\in\mathbb{K}_i$, we have
$k\leftrightsquigarrow k^\prime$ and $p_{kk^{\prime\prime}}^{(n)}=0$ for all $n\ge1$.
\end{itemize}
Then, based on each component $\mathbb{K}_i, 1\le i\le r$, one can define a symbolic system $\theta\colon\varSigma_{\mathbb{K}_i}^+\rightarrow\varSigma_{\mathbb{K}_i}^+$, where $\varSigma_{\mathbb{K}_i}^+$ consists of all $i_{\bcdot}\colon\mathbb{N}\to\mathbb{K}_i$. It is easily seen that $\varSigma_{\mathbb{K}_i}^+$ is a closed invariant subspace of $\varSigma_{\mathbb{K}}^+$.
On the other hand, there holds the following basic result.

\begin{lem}[\cite{Ch60,Sen}]\label{lem2.2}%
Under the basic partition of $\mathbb{K}$ above, there hold the following two statements.
\begin{enumerate}
\item[$\mathrm{(1)}$] $\mu_{\bp,\bP}(\varSigma_{\mathbb{K}_i}^+)>0$, for each $1\le i\le r$.

\item[$\mathrm{(2)}$] $\mathbb{K}_0=\varnothing$ under the assumption $\bp>0$. In general case, $\mathbb{K}_0=\{k\colon p_k=0\}$.
\end{enumerate}
\end{lem}

Let $\alpha_i=\mu_{\bp,\bP}(\varSigma_{\mathbb{K}_i}^+)$ for $1\le i\le r$. Then $0<\alpha_i\le1$ and $\alpha_1+\cdots+\alpha_r=1$. So,
\bean
\varSigma_{\mathbb{K}}^+=\varSigma_{\mathbb{K}_1}^+\cup\cdots\cup\varSigma_{\mathbb{K}_r}^+\quad(\textrm{mod }\mu_{\bp,\bP})
\eean
is a measurable, not necessarily topological, partition of the space $\varSigma_{\mathbb{K}}^+$. Define conditional $\theta$-invariant probability measures $\mu_{\bp,\bP}(\cdot\,\pmb{|}\mathbb{K}_i)$ on $\varSigma_{\mathbb{K}_i}^+$ by
\bean
\mu_{\bp,\bP}(B\,\pmb{|}\mathbb{K}_i)=\frac{\mu_{\bp,\bP}(B\cap\varSigma_{\mathbb{K}_i}^+)}{\alpha_i}\quad\forall B\in\mathscr{B}_{\varSigma_{\mathbb{K}}^+},
\eean
for each $1\le i\le r$. Then,
\bean
\mu_{\bp,\bP}(\cdot)=\alpha_1\mu_{\bp,\bP}(\cdot\,\pmb{|}\mathbb{K}_1)+\cdots+\alpha_r\mu_{\bp,\bP}(\cdot\,\pmb{|}\mathbb{K}_r).
\eean
Next, from Lemma~\ref{lem2.1} one can easily obtain the following standard ergodic decomposition of the $\theta$-invariant probability measure $\mu_{\bp,\bP}$.

\begin{thm}\label{thm2.3}
For each $1\le i\le r$, $\mu_{\bp,\bP}(\cdot\,\pmb{|}\mathbb{K}_i)$ is an ergodic probability measure of $\theta$, restricted on the subspace $\varSigma_{\mathbb{K}_i}^+$.
\end{thm}

Recall that for $\bS=\{S_1,\dotsc,S_K\}\subset\mathbb{R}^{d\times d}$, it is called ``periodically switched stable" (\cite{PR91, Gur95, SWMWK, DHX-aut}) if for any finite-length words $(k_1,\dotsc,k_n)\in\mathbb{K}^n$ and any $n\ge1$, the spectral radius of $S_{k_1}\dotsm S_{k_n}$ is less than $1$, i.e., over any periodical switching sequences
\bean
i_{\bcdot}=(\uwave{k_1,\dotsc,k_n},\uwave{k_1,\dotsc,k_n},\uwave{k_1,\dotsc,k_n},\dotsc)\in\varSigma_K^+,
\eean
we have that $\|S_{i_1}\dotsm S_{i_n}\|\to0$ as $n\to+\infty$.

There are counterexamples which show that the periodically switched stability does not need to imply the absolute asymptotic stability of $\bS$, namely,
$\|S_{i_1}\dotsm S_{i_n}\|\to0$ as $n\to+\infty$ for all $i_{\bcdot}\in\varSigma_K^+$.
See \cite{BM02}, also \cite{BTV03, Koz07, HMST11}, and \cite{DLAA}. However, in \cite{DHX-aut}, the authors proved that $\bS$ is exponentially stable $\mu_{\bp,\bP}$-almost surely, if the transition probability matrix $\bP$ is irreducible, i.e., $\mu_{\bp,\bP}$ is ergodic for the one-sided shift $\theta$. From the ergodic decomposition theorem (Theorem~\ref{thm2.3}) and \cite{DHX-aut}, a more general result can be stated as follows:

\begin{prop}\label{prop2.4}
Let $\bS$ be periodically switched stable. Then the Markovian jump linear system $(\bS,\bxi)$ is exponentially stable $\mathbb{P}$-almost surely; that is to say, for $\mathbb{P}$-a.e. $\omega\in\Omega$,
\bean
\|{S}_{\xi_1(\omega)}\dotsm{S}_{\xi_n(\omega)}\|\xrightarrow[]{\textrm{exponentially fast}}0\quad \textrm{as }n\to+\infty.
\eean
\end{prop}

This generalizes the statement (1) of \cite[Main Theorem]{DHX-aut} from ergodic case to the general case.

\section{Pointwise stabilizability}\label{sec3}%
This section will be devoted to proving our main result.
As in Section~I, we let
\bean
\bS=\{{S}_1,\dotsc,{S}_K\}\subset\mathbb{R}^{d\times d}
\eean
be a set of arbitrarily given $K$ $d$-by-$d$ matrices and
\bean
\bxi=(\xi_n)_{n=1}^{+\infty}\quad \textrm{where }\xi_n\colon(\Omega,\mathscr{F},\mathbb{P})\rightarrow\mathbb{K},
\eean
a discrete-time, time-homogeneous, and stationary $(\bp,\bP)$-Markovian chain, as described before.

\subsection{Main result}
In terms of symbolic dynamics, since $\mu_{\bp,\bP}=\mathbb{P}\circ\varXi^{-1}$, we may state our main result as follows:

\begin{thm}\label{thm3.1}
Let $(\bS,\bxi)$ be a $(\bp,\bP)$-Markovian jump linear system, which is {\it essentially nonuniformly product bounded}, i.e., there exists some function $\beta\colon\Omega\rightarrow[1,\infty)$ such that for $\mathbb{P}$-a.e. $\omega\in\Omega$,
\bean
\|\xi_1(\omega)\dotsm \xi_n(\omega)\|\le\beta(\omega)\quad\forall n\ge1.
\eean
Then the following two statements are equivalent to each other.
\begin{enumerate}
\item[$\mathrm{(a)}^\prime$] $(\bS,\mu_{\bp,\bP})$ is ``pointwise convergent"; that is, to any initial state $x\in\mathbb{R}^{1\times d}$, there corresponds a Borel subset $\pmb{\varSigma}_{x}\subset\varSigma_K^+$ with $\mu_{\bp,\bP}(\pmb{\varSigma}_{x})>0$ such that
\bean
    xS_{i_1}\dotsm S_{i_n}\to\mathbf{0}\;\textrm{ as }\;n\to+\infty,\quad \forall i_{\bcdot}\in\pmb{\varSigma}_{x}.
\eean

\item[$\mathrm{(b)}^\prime$] $(\bS,\mu_{\bp,\bP})$ is ``pointwise exponentially convergent"; i.e., to any initial state $x\in\mathbb{R}^{1\times d}$, there corresponds a Borel subset $\pmb{\varSigma}_{x}^\prime\subset\varSigma_K^+$ with $\mu_{\bp,\bP}(\pmb{\varSigma}_{x}^\prime)>0$ such that
\bean
    \limsup_{n\to+\infty}\frac{1}{n}\log\|xS_{i_1}\dotsm S_{i_n}\|<0,\quad \forall i_{\bcdot}\in\pmb{\varSigma}_{x}^\prime.
\eean
\end{enumerate}
\end{thm}

We will give the proof of Theorem \ref{thm3.1} after we make some remarks. Recall that $\bS$ is itself called ``uniformly product bounded" if the multiplicative semigroup $\bS^+$, generated by $\bS$, is bounded in $\mathbb{R}^{d\times d}$; this is equivalent to that there exists a constant $\beta>0$ such that
\bean
\|\mathrm{S}_{i_1}\cdots\mathrm{S}_{i_n}\|\le\beta\quad\forall n\ge1\textrm{ and }\forall i_{\bcdot}\in\varSigma_K^+.
\eean
It is obvious that the uniform product boundedness of $\bS$ implies the essential nonuniform product boundedness of $(\bS,\mu_{\bp,\bP})$, but the opposite implication does not hold; for example, letting
\bean
\bS=\left\{S_1=\left[\begin{matrix}1&0\\0&1\end{matrix}\right],\; S_2=\left[\begin{matrix}1&1\\0&1\end{matrix}\right]\right\}
\eean
and
$\xi_n(\omega)\in\{1,2\}$, then $(\bS,\mu_{\bp,\bP})$ is essentially nonuniformly product bounded whenever $\bP$ satisfies $p_{12}=0$ and $p_{21}>0$; however, $\bS$ itself is not uniform product bounded. See \cite{DHX-pre} for more stricter counterexamples.

A vector norm $|\pmb{|}\cdot\pmb{|}|_*$ on $\mathbb{R}^{1\times d}$ is called a ``pre-extremal" norm of $\bS$, if its induced matrix norm on $\mathbb{R}^{d\times d}$ is such that $|\pmb{|}\mathrm{S}_i\pmb{|}|_*\le1$ for all $i\in\{1,\dotsc,K\}$. If $\bS$ is uniform product bounded, then such a pre-extremal norm always exists; see, for example, \cite{BT79, Koz90, El95, OR97, Hart02, Dai-JDE}. However, in the situation of Theorem~\ref{thm3.1}, there does not have a pre-extremal norm for $\bS$ in general \cite{DHX-pre}.
In addition, since here the dimension $d$ of $\mathbb{R}^{1\times d}$ is not less than $2$, for any initial state $x\in\mathbb{R}^{1\times d}$ with $\|x\|=1$ and any $i_{\bcdot}\in\varSigma_K^+$, as $n\to+\infty$
\bean
    xS_{i_1}\dotsm S_{i_n}\to\mathbf{0}\quad\not\Rightarrow\quad S_{i_1}\dotsm S_{i_n}\to\mathbf{0};
\eean
and crucially the sequence $f_n(i_{\bcdot}):=\log\|xS_{i_1}\dotsm S_{i_n}\|$ of functions of the variable $i_{\bcdot}\in\varSigma_K^+$ does not have the subadditivity property. So, we cannot apply the Subadditive Ergodic Theorem to the proof of Theorem~\ref{thm3.1} in the general case.

In the recent paper~\cite{DHX-pre}, using a splitting of the state space we have proven the following.

\begin{prop}[\cite{DHX-pre}]\label{prop3.3}
Let $(\bS,\bxi)$ be a $(\bp,\bP)$-Markovian jump linear system. If $(\bS,\bxi)$ is essentially nonuniformly product bounded the following two statements are equivalent to each other.
\begin{itemize}
\item $(\bS,\bxi)$ is ``consistently convergent"; that is, there is measurable set $\Omega^\prime\subset\Omega$ with $\mathbb{P}(\Omega^\prime)>0$ such that
\bean
{S}_{\xi_1(\omega)}\dotsm S_{\xi_n(\omega)}\to\mathbf{0}\; \textrm{as }n\to+\infty,\quad \forall \omega\in\Omega^\prime.
\eean

\item $(\bS,\bxi)$ is ``consistently exponentially convergent"; that is, there exists a measurable set $\Omega^{\prime\prime}\subset\Omega$ with $\mathbb{P}(\Omega^{\prime\prime})>0$ such that as $n\to+\infty$,
\bean
{S}_{\xi_1(\omega)}\dotsm S_{\xi_n(\omega)}\xrightarrow[]{\textrm{exponentially fast}}\mathbf{0}\quad\forall \omega\in\Omega^{\prime\prime}.
\eean
\end{itemize}
\end{prop}

\begin{proof}
If we assume the more stronger condition that $\bS$ is uniformly product bounded, then using Kingman's subadditive ergodic theorem in ergodic theory and Egoroff's almost uniform convergence theorem in measure theory, we can simply prove this statement as follows.

According to Theorem~\ref{thm2.3}, there is no loss of generality in assuming that $\mu_{\bp,\bP}$ is ergodic for the Markovian shift $\theta\colon i_{\bcdot}\mapsto i_{\bcdot+1}$ on $\varSigma_K^+$. Let there be a Borel subset $\pmb{\varSigma}^\prime$ of $\varSigma_K^+$ with $\mu_{\bp,\bP}(\pmb{\varSigma}^\prime)>0$ such that
\bean
\|S_{i_1}\dotsm S_{i_n}\|\to0\;\textrm{ as }n\to+\infty,\quad \textrm{for }\mu_{\bp,\bP}\textrm{-a.e. }i_{\bcdot}\in\pmb{\varSigma}^\prime.
\eean
Since $(\bS,\bxi)$ is $\mu_{\bp,\bP}$-essentially bounded, there is a $\beta>1$ such that for $\mu_{\bp,\bP}$-a.e. $i_{\bcdot}\in\varSigma_K^+$,
\bean
\|S_{i_1}\dotsm S_{i_n}\|\le\beta\quad\forall n\ge1.
\eean
Take a number $0<\alpha<1$ which is so small that
\bean
\frac{1}{2}\mu_{\bp,\bP}(\pmb{\varSigma}^{\prime})\log\alpha+\log\beta<0.
\eean
From Egoroff's almost uniform convergence theorem, we can take a Borel set $\pmb{\varSigma}^{\prime*}\subseteq\pmb{\varSigma}^{\prime}$ with $\mu_{\bp,\bP}(\pmb{\varSigma}^{\prime*})>\frac{1}{2}\mu_{\bp,\bP}(\pmb{\varSigma}^{\prime})$ and an integer $N\ge1$ such that
\bean
\|S_{i_1}\dotsm S_{i_n}\|\le\alpha\quad\forall n\ge N\textrm{ and }i_{\bcdot}\in\pmb{\varSigma}^{\prime*}.
\eean
Then from the subadditive ergodic theorem, see e.g. \cite[Theorem~10.1]{Wal82}, it follows that for $\mu_{\bp,\bP}$-a.e. $i_{\bcdot}\in\varSigma_K^+$ and any $n\ge N$,
\begin{align*}
&\lim_{m\to+\infty}\frac{1}{m}\log\|S_{i_1}\dotsm S_{i_m}\|\\
&\qquad=\inf_{m\ge1}\frac{1}{m}\int_{\varSigma_K^+}\log\|S_{i_1}\dotsm S_{i_m}\|\mu_{\bp,\bP}(di_{\bcdot})\\
&\qquad\le\frac{1}{n}\int_{\varSigma_K^+}\log\|S_{i_1}\dotsm S_{i_n}\|\mu_{\bp,\bP}(di_{\bcdot})\\
&\qquad\le\frac{1}{n}\left(\frac{1}{2}\mu_{\bp,\bP}(\pmb{\varSigma}^{\prime})\log\alpha+\log\beta\right)\\
&\qquad<0.
\end{align*}
This completes the proof of Proposition~\ref{prop3.3}.
\end{proof}

Recall, for instance from \cite{SU94} and \cite{Sun04} in the deterministic situation, that $\bS$ is called to be:
\begin{itemize}
\item ``consistently convergent" if there is a switching sequence $i_{\bcdot}\in\varSigma_K^+$ such that
\bean
 S_{i_1}\dotsm S_{i_n}\to\mathbf{0}\; \textrm{as }n\to+\infty,
\eean
or equivalently,
\bean
xS_{i_1}\dotsm S_{i_n}\to\mathbf{0}\;  \textrm{as }n\to+\infty,\quad\forall x\in\mathbb{R}^{1\times d};
\eean

\item ``consistently exponentially convergent" if there exists a switching sequence $i_{\bcdot}\in\varSigma_K^+$ such that
\bean
\limsup_{n\to+\infty}\frac{1}{n}\log\|{S}_{i_1}\dotsm S_{i_n}\|<0\; \textrm{as }n\to+\infty,
\eean
that is to say,
\bean
{S}_{i_1}\dotsm S_{i_n}\xrightarrow[]{\textrm{exponentially fast}}\mathbf{0}\; \textrm{as }n\to+\infty.
\eean
\end{itemize}
In D.\,P.~Stanford and J.\,M.~Urbano~\cite[Theorem~3.5]{SU94}, it was proved that $\bS$ is consistently convergent if and only if it is consistently exponentially convergent; more precisely, $\bS$ is consistently convergent if and only if there is a finite-length word $w=(k_1,\dotsc,k_m)$ in $\mathbb{K}^m$, for some $m\ge1$, such that the spectral radius $\rho({S}_{k_1}\cdots{S}_{k_m})<1$. Also see Z.~Sun \cite[Proposition~4]{Sun04} and J.-W.~Lee and G.\,E.~Dullerud \cite[Theorem~2]{LD07}.

We notice here that although the consistent exponential convergence of $\bS$ implies, from Y.~Huang \textit{et al.} \cite{HLHM}, that there exists some other $(\bp^\prime, \bP^\prime)$-Markovian probability measure $\mu_{\bp^\prime, \bP^\prime}$ such that $\bS$ is exponentially convergent $\mu_{\bp^\prime, \bP^\prime}$-almost surely, yet it cannot imply the consistent exponential convergence of $(\bS,\bxi)$ in general. This is because $\mu_{\bp^\prime, \bP^\prime}$, constructed in \cite{HLHM} there, does not need to equal $\mu_{\bp, \bP}$ that has been presented by $\bxi$ in our situation, and the set of all periodical switching sequences in $\varSigma_K^+$ has $\mu_{\bp,\bP}$-measure $0$ in general case; see for example, \cite{DHX-aut}.
\subsection{Proof of Theorem~\ref{thm3.1}}\label{sec3.2}
For any nonempty subspace $E$ of $\mathbb{R}^{1\times d}$ and any $S_i\in\bS$, we write $S_i(E)=\{xS_i\,|\,x\in E\}$. To prove Theorem~\ref{thm3.1} we need the following dichotomy theorem, which comes directly from \cite{DHX-pre}.

\begin{lem}[\cite{DHX-pre}]\label{lem3.4}
Let $(\bS,\bxi)$ be a $(\bp,\bP)$-Markovian jump linear system, which is essentially nonuniformly product bounded.
Then, there exists an $\theta$-invariant Borel subset $\pmb{\varSigma}_{\bp,\bP}^+$ of $\varSigma_K^+$ with $\mu_{\bp,\bP}$-measure $1$ such that for any
$i_{\bcdot}\in\pmb{\varSigma}_{\bp,\bP}^+$, there corresponds a direct sum decomposition of $\mathbb{R}^{1\times d}$ into subspaces
\bean
\mathbb{R}^{1\times d}=E^s(i_{\bcdot})\oplus E^c(i_{\bcdot}),
\eean
where $i_{\bcdot}\mapsto E^s(i_{\bcdot})$ is Borel measurable, with the invariance $S_{i_1}(E^s(i_{\bcdot}))\subseteq E^s(i_{\bcdot+1})$ for all $i_{\bcdot}\in\pmb{\varSigma}_{\bp,\bP}^+$, for which there hold the following two properties:
\begin{enumerate}
   \item[$(1)$] for any $i_{\bcdot}\in\pmb{\varSigma}_{\bp,\bP}^+$, one can find a positive integer sequence $n_k(i_{\bcdot})\nearrow+\infty$ with $i_{\bcdot+n_k}\to i_{\bcdot}$ such that
\bean
S_{i_1}\dotsm S_{i_{n_k}}\!\upharpoonright\!E^c(i_{\bcdot})\to\mathrm{Id}_{\mathbb{R}^{1\times d}}\!\upharpoonright\!E^c(i_{\bcdot})\quad \textrm{as }k\to+\infty;
\eean

\item[$(2)$] for any $i_{\bcdot}\in\pmb{\varSigma}_{\bp,\bP}^+$,
\bean
\lim_{n\to+\infty}\frac{1}{n}\log\|x S_{i_1}\dotsm S_{i_n}\|<0\quad\forall x\in E^s(i_{\bcdot})
\eean
and
\bean
\limsup_{n\to+\infty}\|x S_{i_1}\dotsm S_{i_n}\|>0\quad \forall x\in \mathbb{R}^{1\times d}\setminus E^s(i_{\bcdot}).
\eean
\end{enumerate}
\end{lem}
Its proof involves ergodic theory and semigroup theory.

We make some comments here. Firstly, we cannot guarantee the invariance of the central manifold $E^c(i_{\bcdot})$, since the Markovian shift $\theta\colon i_{\bcdot}\mapsto i_{\bcdot+1}$ is not homeomorphic here.

Secondly, Lemma~\ref{lem3.4}.(2) implies that
\bean
\|S_{i_1}\dotsm S_{i_n}\!\upharpoonright\!E^s(i_{\bcdot})\|\to0\;\textrm{ as }\;n\to+\infty\quad\forall i_{\bcdot}\in\pmb{\varSigma}_{\bp,\bP}^+;
\eean
and for any $\varepsilon>0$, there is a Borel subset $\pmb{\varSigma}_{\bp,\bP}^{+\prime}\subset\pmb{\varSigma}_{\bp,\bP}^+$ with $\mu_{\bp,\bP}(\pmb{\varSigma}_{\bp,\bP}^{+\prime})>1-\varepsilon$ such that
\bean
\|S_{i_1}\dotsm S_{i_n}\!\upharpoonright\!E^s(i_{\bcdot})\|\to0
\eean
uniformly for $i_{\bcdot}\in\pmb{\varSigma}_{\bp,\bP}^{+\prime}$ from Egoroff's almost uniform convergence theorem.
However, in general, we cannot expect $\mu_{\bp,\bP}(\pmb{\varSigma}_{\bp,\bP}^{+\prime})=1$. This is just the essential difference between ``uniformity" and ``non-uniformity" in the smooth ergodic theory and linear cocycle theory.

Thirdly, if $(\bS,\bxi)$ is pointwise convergent, then for any nonzero $x\in\mathbb{R}^{1\times d}$ we have form Lemma~\ref{lem3.4} that
\bean
\mu_{\bp,\bP}(\pmb{\varSigma}_{x}^+)>0\quad\textrm{where }\pmb{\varSigma}_{x}^+:=\left\{i_{\bcdot}\in\pmb{\varSigma}_{\bp,\bP}^+\colon x\in E^s(i_{\bcdot})\right\}.
\eean
However, as in Remark~\ref{rem1.3},
$\mu_{\bp,\bP}(\pmb{\varSigma}_{x}^+)\not=1$
in general case, even though $\mu_{\bp,\bP}$ is ergodic; since there is no the invariance that $\theta(\pmb{\varSigma}_{x}^+)\subseteq\pmb{\varSigma}_{x}^+$, unlike the stable set $\Omega_{\mathrm{stable}}$ defined in Remark~\ref{rem1.3}, unless $\bS$ is of diagonal form.

The pointwise exponential convergence of $(\bS,\bxi)$ implies obviously the pointwise convergence. Thus, according to the ergodic decomposition (Theorem~\ref{thm2.3}), Theorem~\ref{thm3.1} follows immediately from the following statement.

\begin{lem}\label{lem3.5}%
Let $(\bS,\bxi)$ be a $(\bp,\bP)$-Markovian jump linear system, which is essentially nonuniformly product bounded and ergodic. If $(\bS,\mu_{\bp,\bP})$ is pointwise convergent, then it is pointwise exponentially convergent.
\end{lem}

\begin{proof}
This statement comes at once from Lemma~\ref{lem3.4}.
In fact, let $x\in\mathbb{R}^{1\times d}\setminus\{\mathbf{0}\}$ be arbitrary. Since $(\bS,\mu_{\bp,\bP})$ is pointwise convergent, one can find some Borel subset $\pmb{\varSigma}_{x}$ of $\varSigma_K^+$
with $\mu_{\bp,\bP}(\pmb{\varSigma}_{x})>0$ such that
\bean
\|x{S}_{i_1}\cdots{S}_{i_n}\|\to0\;\textrm{ as }n\to+\infty\quad \forall i_{\bcdot}\in\pmb{\varSigma}_{x}.
\eean
According to Lemma~\ref{lem3.4}, one can further choose a Borel subset $\pmb{\varSigma}_{x}^\prime\subseteq\pmb{\varSigma}_{x}$ with $\mu_{\bp,\bP}(\pmb{\varSigma}_{x}^\prime)>0$ such that $x\in E^s(i_{\bcdot})$ for any $i_{\bcdot}\in\pmb{\varSigma}_{x}^\prime$. This leads to the conclusion.
\end{proof}

Therefore the proof of Theorem~\ref{thm3.1} is completed.
\section{Concluding remarks}\label{sec4}
For a discrete-time Markovian jump linear system, in this note we have introduced two concepts --- pointwise convergence and pointwise exponential convergence in terms of probability measure. These two types of convergences in general are not equivalent to each other. The classification of the type of convergences is important in many aspects, such as in numerical computations, in optimal control and so on. Thus it it important to know under what condition they share the same convergent property. However, to show that is not so straightforward. In this note, we proved that if the Markovian jump linear system is essentially nonuniformly product bounded, then the pointwise convergence and pointwise exponential convergence are equivalent to each other under the framework of symbolic dynamics.

\section*{Acknowledgment}
The authors would like to thank the anonymous referees for their critical comments and suggestions that lead to the improvement of this manuscript.

This work was supported partly by National Natural Science Foundation of China (Grant Nos.~11071112, 11071263), NSF of Guangdong Province and in part by NSF 1021203 of the United States.

\ifCLASSOPTIONcaptionsoff
  \newpage
\fi

\end{document}